\documentclass[12pt]{article}
\usepackage{amsmath}
\usepackage{graphicx,psfrag,epsf}
\usepackage{enumerate}
\usepackage{url} 

 \newtheorem{theo}{\small\bf Theorem}
 \newtheorem{lem}{\small\bf Lemma}
 \newtheorem{rem}{\small\bf Remark}
 \newenvironment{REM}{\begin{rem} \rm}{\end{rem}}
 \newtheorem{exam}{\small\bf Example}
 
 \newtheorem{defi}{\small\bf Definition}
 \newenvironment{DEFI}{\begin{defi} \rm}{\end{defi}}
 \newtheorem{cor}{\small\bf Corollary}
 \renewcommand{\Pr}{\mbox{\rm  \hspace*{.2ex}I\hspace{-.5ex}P\hspace*{.2ex}}}
 \newcommand{\be}{\begin{equation}}
 \newcommand{\ee}{\end{equation}}
 \newcommand{\R}{\mbox{R}}
 \newcommand{\Z}{\mbox{Z}}
 \newcommand{\E}{\mbox{\rm \hspace*{.2ex}I\hspace{-.5ex}E\hspace*{.2ex}}}

 \newenvironment{pr}[1]{{\small\bf {#1}:}}{}

 \renewcommand{\Box}{\mbox{Q.E.D.}}

   \font\twelvebm                       = cmmib10 at 12truept
   \font\tenbm                          = cmmib10 at 10truept
   \font\sevenbm                        = cmmib10 at 7truept

\tenbm \scriptscriptfont9              = \sevenbm

  at 10truept
 at 10truept  at 10truept  at 10truept 
at 10truept

\begin{document}

\def\spacingset#1{\renewcommand{\baselinestretch}%
{#1}\small\normalsize} \spacingset{1}



 \title{\bf \noindent

 \centerline{\LARGE\bf The characteristic function}

  \centerline{\LARGE\bf of the discrete Cauchy distribution}
  \bigskip

  \centerline{\it In Memory of T.\ Cacoullos}
  }
  \bigskip

 \author{{\Large Nickos Papadatos}
    \hspace{.2cm}
    \\
    Department of Mathematics, National and Kapodistrian
    University of Athens,
    \\
    Panepistemiopolis, 157 84 Athens, Greece}

 \date{}





  \maketitle

 \thispagestyle{empty}

\bigskip
\begin{abstract}
 A new family of integer-valued Cauchy-type
 distributions is introduced,
 the {\it Cauchy-Cacoullos family}.
 The characteristic function is evaluated,
 showing some interesting distributional
 properties, similar to the ordinary (continuous) Cauchy scale family.
 The results are extendable to discrete Student-type distributions
 with odd degrees of freedom.
\end{abstract}

\noindent%
\vfill

\spacingset{1} 
\section{Introduction and summary}
\label{sec.1}

 Some years ago, Cacoullos (Personal Communication),
 considering discretization of well-known continuous distributions,
 introduced a (standard)
 discrete Cauchy random variable (r.v.)\ $X$
 with probability mass function
 (p.m.f.)
 \be
 \label{1}
 \Pr(X=k)=\frac{1/\pi_0}{1+k^2}, \ \ \ k\in\Z,
 \ee
 by the obvious substitution $k\in\Z$ for $x\in\R$
 in the standard Cauchy density
 \be
 \label{1b}
 f(x)=\frac{1/\pi}{1+x^2}, \ \ x\in\R.
 \ee

 Cacoullos immediately raised two natural questions:

 \begin{itemize}
 \item[(A)] While it is expected to be very close to $\pi$,
 what is the exact value of the normalizing constant
 $\pi_0$ in (\ref{1})?
 \smallskip

 \item[(B)] While the characteristic function (ch.f.)\ of (\ref{1b})
 is $\phi(t)=e^{-|t|}$,
 what is the corresponding one, say $\phi_1$,
 of (\ref{1})?
 \end{itemize}

 We provide explicit answers in section \ref{sec.2}.
 It is well-known that the
 (continuous) Cauchy
 distribution appears naturally in statistics and probability.
 At this point it should be noted though the standard Cauchy
 r.v.\ is customarily defined
 as the ratio of two independent standard normal r.v.'s,
 or as the tangent of a randomly chosen angle in $[0,2\pi)$,
 it has recently been shown (\cite{1}, \cite{8},
 \cite{10}) that the ratio representation still holds if $(X,Y)$ follows
 any bivariate spherically
 symmetric distribution.

 In \cite{4}, Cacoullos showed that if $X=(X_1,\ldots,X_p)'$
 ($p\geq 3$) is spherically symmetrically distributed around zero
 then all polar angle tangent vectors follow a multivariate
 Cauchy; note that, e.g., Feller (1966) defines the symmetric
 bivariate and trivariate Cauchy distributions directly
 through their densities -- not as tangent vectors.




 In contrast to (\ref{1b}) and its location-scale
 extension, for which several applications are known
 both in probability and statistics, for (\ref{1}) we
 have been able to find few results related to
 stochastic processes -- see, e.g., \cite{11}, p.\ 383.
 However, the asymptotic distribution of the sample means
 for (\ref{1}), Theorem \ref{theo.2}, may serve as a
 starting point for applications; so appears to be
 the Cauchy-Cacoullos family defined by (\ref{3}).
 These considerations are, however, beyond the scope of the
 present note.



 In section \ref{sec.3}
 we introduce a novel family of integer-valued distributions,
 the
 {\it Cauchy-Cacoullos family},
 sharing similar 
 properties --
 see Definition \ref{def.1} and Remark \ref{rem.3}.
 In particular, any distribution in this family has a simple
 characteristic function that can be written down explicitly,
 Theorem \ref{theo.3}, and the same is valid for the discrete Student-type
 distributions of Remark \ref{rem.3}.
 Basic inference properties for this family are included in
 Theorem \ref{theo.3b}, while some distributional properties
 are discussed in some detail in Section \ref{sec.4}; see Theorems \ref{theo.2}--\ref{theo.6}.
 We hope that the proposed simple formulae will enlarge
 the applicability of discrete Cauchy distribution in the future.

 \section{The characteristic function}
 \label{sec.2}

 Since $\phi_1(t)=\E e^{itX}=\E\cos(tX)+i\E\sin(tX)$ ($i$ denotes the imaginary
 unit) and $X$
 is symmetrically distributed around the origin (hence, $\E\sin(tX)=0$), both
 questions, (A), (B), will be answered if we manage to calculate in
 a closed form the function
 $g:\R\to\R$,
 defined by the Fourier series
 \be
 \label{2}
 g(t):=\sum_{n=0}^{\infty} \frac{\cos(n t)}{1+n^2},
 \  \ \ \
 t\in\R.
 \ee

 Therefore, the problem is to identify which function $g$ is represented
 as a series of cosines
 with Fourier coefficients as in (\ref{2}). Clearly, $g$ is periodic with period
 $2\pi$. Thus, it suffices to restrict
 our attention to $t$-values in
 the interval $-\pi\leq t\leq \pi$. On the other hand, since a
 cosine Fourier series corresponds to an even function, we
 may further restrict the $t$-values into the interval
 $0\leq t\leq \pi$.

 The key lemma is:

 \begin{lem}
 \label{lem.1}
 For $-2\pi\leq t\leq 2\pi$,
 \[
 g(t)=\frac{1}{2}+\frac{\pi\cosh(\pi-|t|)}{2\sinh(\pi)}.
 \]
 \end{lem}

 \noindent
 We omit the proof because we shall show a more general result in
 Section \ref{sec.3}, below.

 \begin{cor}
 \label{cor.1}
 The normalizing constant $\pi_0$ is given by
 \[
 \pi_0=2g(0)-1=\frac{\pi\cosh(\pi)}{\sinh(\pi)}
 =\pi\left(1+\frac{2}{e^{2\pi}-1}\right)\simeq3.15334809493716\ldots\ .
 \]
 \end{cor}

\begin{figure}[htp]
 \centering
 \includegraphics[width=.6\linewidth]{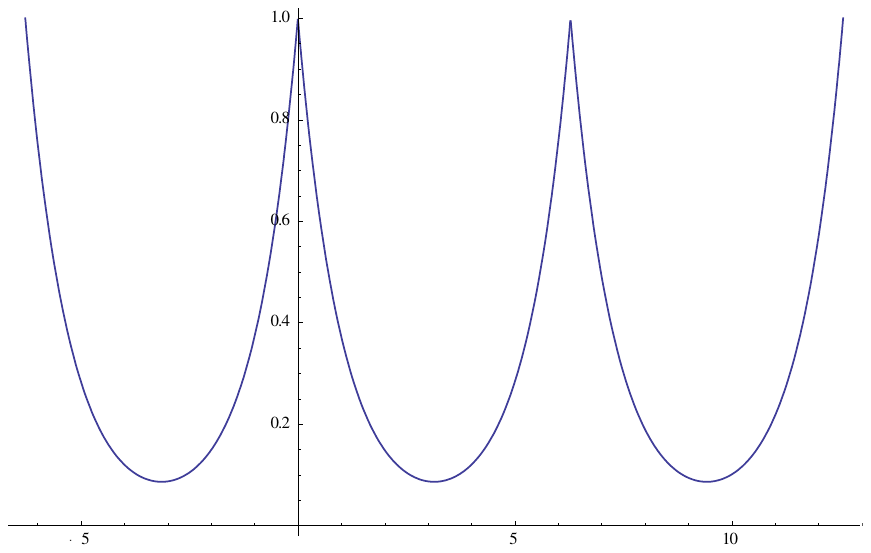}
 \vspace*{-46ex}
 \caption{The characteristic function $\phi_1(t)$
 in the interval $-2\pi\leq t\leq 4\pi$.}
 \label{fig.1}
 \end{figure}

 The formula for the ch.f., and is an immediate consequence of Lemma
 \ref{lem.1} and (\ref{2}):
 \begin{theo}
 \label{theo.1}
 The ch.f.\ of $X$ is given by
 $
 \phi_1(t)=\cosh(\pi-|t|)/\cosh(\pi)$,
 $-2\pi\leq t \leq 2\pi$,
 and it is periodic with period $2\pi$ {\rm (see Fig.\ \ref{fig.1})}.
 \end{theo}



 \vspace*{-3em}


 \section{The Cauchy-Cacoullos family of discrete distributions}
 \label{sec.3}
  If we multiply a continuous Cauchy r.v.\ by a constant $\lambda>0$
 we stay in the same family of distributions -- the Cauchy scale family.
 More precisely, if $X$ is standard Cauchy, the density of $\lambda X$
 is given by
 \[
 f(x)= \frac{1}{\pi}\ \frac{\lambda}{\lambda^2+x^2},  \ \ \ x\in\R, \ \ \lambda>0.
 \]
 However, this is no longer true for a discrete Cauchy $X$, since the support
 of $\lambda X$ is not the set of integers.
 Motivated from this observation, we
 define a family of discrete integer-valued distributions as
 follows:
 \begin{DEFI}
 \label{def.1}
 The discrete Cauchy-Cacoullos family (${\cal CC}$, for short)
 contains the p.m.f.'s
 \be
 \label{3}
 f_{\lambda}(k)=\frac{\tanh(\lambda\pi)}{\pi}
 \ \frac{\lambda}{\lambda^2+k^2},
 \ \ \ \ k\in\Z, \ \ \lambda>0.
 \ee
 For completeness of the presentation, it is convenient to
 include the limiting case $\lambda=0$, which corresponds
 to a degenerate r.v.\ at zero.
 \bigskip
 \end{DEFI}

 \noindent
 Although this family has several interesting properties, similar
 to the Cauchy, it does not seem to have been studied
 elsewhere. Clearly, for $\lambda=1$ we get (\ref{1}). At a first
 glance, it is
 not entirely obvious to verify that the normalizing constant is as in
 (\ref{3}).
 This is a by-product of the following result.
 \begin{lem}
 \label{lem.2}
 For $-\pi\leq t\leq \pi$ and $\lambda>0$,
 \[
 \cosh(\lambda t)=\frac{\lambda \sinh(\lambda\pi)}{\pi}
 \left\{\frac{1}{\lambda^2}+2\sum_{n=1}^\infty \frac{(-1)^n\cos(nt)}{\lambda^2+n^2}
 \right\}.
 \]
 \end{lem}
 \begin{pr}{Proof} We express the even function $h(t)=\cosh(\lambda t)$
 in a cosine Fourier series to get
 $h(t)\sim \sum_{n=0}^\infty \alpha_n \cos(nt)$.
 Simple calculations show that
 \[
 \alpha_0=\frac{1}{2\pi}\int_{-\pi}^\pi \cosh(\lambda u)du=
 \frac{\sinh(\lambda\pi)}{\pi\lambda}
 \]
 and
 \[
 \alpha_n=\frac{1}{\pi}\int_{-\pi}^\pi \cosh(\lambda u)\cos(n u)du
 =
 (-1)^n
 \frac{2\lambda\sinh(\lambda\pi)}{\pi(\lambda^2+n^2)}, \ \
 n=1,2,\ldots \ .
 \]
 Since $h$ is differentiable in $[-\pi,\pi]$ with $h(-\pi)=h(\pi)$,
 the lemma is proved (and the series converges uniformly to $h$).
 \bigskip
 $\Box$
 \end{pr}

 \noindent
 If we set $\lambda=1$ and $t\to t-\pi$ in Lemma \ref{lem.2} we obtain Lemma \ref{lem.1}
 with $g$ as in (\ref{2}).

  \begin{cor}
 \label{cor.2} We have
 \[
 \sum_{k=-\infty}^\infty \frac{1}{\lambda^2+k^2}
 =\frac{\pi}{\lambda\tanh(\lambda \pi)},
 \]
 and hence, {\rm (\ref{3})} defines a p.m.f.\ for any $\lambda>0$.
 \end{cor}
 \begin{pr}{Proof}
 Substitute $t=\pi$ in Lemma {\ref{lem.2}}.
 $\Box$
 \bigskip
 \end{pr}

 As for the case $\lambda=1$, we can obtain the ch.f.\ of
 $X_\lambda\sim f_{\lambda}$ in a closed form.
 \begin{theo}
 \label{theo.3}
 The ch.f.\ of $X_\lambda$ with p.m.f.\ $f_\lambda\in \mbox{\rm \mbox{${\cal CC}$}}$ is given by
 \[
 \phi_\lambda (t)=\frac{\cosh\big(\lambda(t-\pi)\big)}{\cosh(\lambda\pi)}, \ \ \
 0\leq t \leq 2\pi,
 \]
 and it is periodic with period $2\pi$. More precisely,
 \[
 \phi_\lambda(t)=\frac{\cosh\left(\lambda\Big(t-2\pi\lfloor
 \frac{t}{2\pi}\rfloor-\pi\Big)\right)}{\cosh(\lambda\pi)},
 \ \ -\infty<t<\infty,
 \]
 where $\lfloor x\rfloor$ denotes the integer part of $x$.
 \end{theo}
 \begin{pr}{Proof} As is well-known,
 all integer-valued r.v.'s have periodic
 ch.f.'s, with period $2\pi$. The particular r.v.\ is
 symmetrically distributed around zero, and thus, its ch.f.\ is real and even,
 so that $\phi_\lambda(t)=\E\cos(t X_\lambda)$.
 To calculate this, we may restrict our
 attention in the interval $0\leq t\leq 2\pi$. Then, since  $-\pi\leq t-\pi\leq \pi$
 and $\cos(nt)=(-1)^n \cos(n (t-\pi))$,
 \begin{eqnarray*}
 \phi_{\lambda}(t)
 & =&
 \frac{\lambda\tanh(\lambda\pi)}{\pi}
 \left\{\frac{1}{\lambda^2}+2 \sum_{n=1}^\infty
 \frac{(-1)^n \cos(n (t-\pi))}{\lambda^2+n^2}\right\}
 \\
 &=&
 \frac{\lambda\tanh(\lambda\pi)}{\pi}
 \
 \frac{\pi\cosh(\lambda(t-\pi))}{\lambda\sinh(\lambda\pi)},
 \end{eqnarray*}
 where the second equality follows from Lemma \ref{lem.2}.
 $\Box$
 \bigskip
 \end{pr}

 Statistical inference for the parameter $\lambda$
 is facilitated  from the fact that the p.m.f.'s and the ch.f.'s in \mbox{${\cal CC}$}
 have tractable forms.

 \begin{theo}
 \label{theo.3b}
 Consider a random sample $X_1,\ldots,X_n\sim f_\lambda\in {\cal CC}$
 with $\lambda>0$ unknown.
 \\
 {\rm (i)} The minimal sufficient statistic is $T=(Y_1,\ldots,Y_n)$,
 with $Y_1\leq Y_2\leq \cdots\leq Y_n$ being the order statistics
 of $|X_1|,\ldots,|X_n|$. \\
 {\rm (ii)} The Fisher Information {\rm (of a single observation)} is
 \be
 \label{I}
 I(\lambda)=\frac{1}{2\lambda^2}+\frac{\pi}{\lambda} w(\lambda)
 \ \ \mbox{where} \ \ w(\lambda)=\frac{\lambda\pi}{\cosh(\lambda\pi)^2}
 -\frac{1}{\sinh(2\lambda\pi)}.
 \ee
 {\rm (iii)} The MLE $\widehat{\lambda}_n$ of $\lambda$ is unique; it is given as
 the unique solution in $[0,\infty)$ of the equation
 \be
 \label{MLE}
 \frac{\pi\lambda}{\sinh(2\pi\lambda)}+\frac{1}{n}\sum_{i=1}^n\frac{X_i^2}{\lambda^2+X_i^2}=\frac{1}{2}.
 \ee
 {\rm (iv)} The MLE is consistent and asymptotically efficient,
 \[
 \sqrt{n}\left(\widehat{\lambda}_n-\lambda\right)\stackrel{d}{\to} N(0,1/I(\lambda)),
 \]
 where $\stackrel{d}{\to}$ denotes weak convergence.
 \end{theo}
 \begin{pr}{Proof}
 Let ${\bf x}=(x_1,\ldots,x_n)$ and ${\bf y}=(y_1,\ldots,y_n)$ be two vectors
 in $\Z^n$. Then, the likelihood ratio is given by
 \[
 \frac{L({\bf x};\lambda)}{L({\bf y};\lambda)}=\prod_{i=1}^n \frac{\lambda^2+y_i^2}{\lambda^2+x_i^2},
 \]
 and it has the same form  as in the continuous Cauchy scale-family. Obviously, this ratio is independent of
 $\lambda>0$ if and only if the ordered squared values of ${\bf x}$ and ${\bf y}$ are identical, and this verifies (i). Now, a straightforward computation yields the score function
 \[
 S(k;\lambda):=\frac{\partial}{\partial \lambda} \log f_\lambda(k)=\frac{1}{\lambda}
 +\frac{2\pi}{\sinh(2\pi\lambda)}-\frac{2\lambda}{\lambda^2+k^2}.
 \]
 Let $X_\lambda\sim f_\lambda$.
 Using Remarks \ref{rem.2}, \ref{rem.3} below, it is seen that $\E S(X_\lambda;\lambda)=0$ and
 $\E S(X_\lambda;\lambda)^2=I(\lambda)$ with $I(\lambda)$
 as in (\ref{I}). Note that $I(\lambda)=-\E\frac{\partial^2}{\partial \lambda^2} \log f_\lambda(X_\lambda)$, since the regularity conditions
 are obviously fulfilled; both formulae require computation of the series $\sum_n (\lambda^2+n^2)^{-s}$, $s=1,2,3$. Moreover, one can easily verify that the log-likelihood is given by
 \be
 \label{Like}
 \frac{\partial}{\partial\lambda}\log L({\bf x};\lambda)=\frac{2n}{\lambda}
 \left(\frac{\pi\lambda}{\sinh(2\pi\lambda)}-\frac{1}{2}+\frac{1}{n}\sum_{i=1}^n \frac{x_i^2}{\lambda^2+x_i^2}\right).
 \ee
 For fixed ${\bf x}\in \Z^n$, the positive function $u(\lambda):=\pi\lambda/\sinh(2\pi\lambda)+n^{-1}
 \sum_{i=1}^n x_i^2/(\lambda^2+x_i^2)$ decreases to zero as $\lambda\to\infty$ and has a limit
  $u(0+)\geq 1/2$ (it equals to $1/2$ iff ${\bf x}={\bf 0}$).
  Since $u$ is strictly decreasing and continuous, the likelihood is first increasing and
  then decreasing, reaching its global maximum at $\lambda_0$, where $u(\lambda_0)=1/2$.
  This shows that the MLE is the unique solution of (\ref{MLE}),
  it equals to $0$ iff ${\bf X}={\bf 0}$, and it is otherwise positive.
  Finally, in order to prove (iv), fix $\lambda=\lambda_0$ and $c\in(0,\lambda_0)$, and
  assume that $\lambda$ varies in the interval $(\lambda_0-c,\lambda_0+c)$.
  Then, $\frac{\partial^3}{\partial\lambda^3}\log f_{\lambda}(k)=A(\lambda)+B(\lambda,k)$
  where
  \[
  A(\lambda)= 4\pi^3\ \frac{3+\cosh(4\lambda\pi)}{\sinh(2\lambda\pi)^3}+\frac{2}{\lambda^3},
  \ \
  B(\lambda,k)=4\lambda \frac{3k^2-\lambda^2}{(\lambda^2+k^2)^3}.
  \]
  The function $A$ is decreasing and positive, so that $|A(\lambda)|<A(\lambda_0-c)$.
  Moreover,
  \[
  |B(\lambda,k)|<4\lambda \frac{3k^2+3\lambda^2}{(\lambda^2+k^2)^3}
  <\frac{12(\lambda_0+c)}{((\lambda_0-c)^2+k^2)^2}\leq \frac{12(\lambda_0+c)}{(\lambda_0-c)^4}.
  \]
  It follows that we can find a finite constant $M=M(\lambda_0,c)$
  such that $|\frac{\partial^3}{\partial\lambda^3}\log f_{\lambda}(k)|<M$
  uniformly in $k\in\Z$, $\lambda\in(\lambda_0-c,\lambda_0+c)$, and the result
  follows by applying Theorem 3.10 in \cite{LG-N9}.
 $\Box$
 \bigskip
 \end{pr}

 Unfortunately, the MLE does not admit a closed form and, hence, numerical procedures
 should be employed. On the other hand, we can construct closed-form consistent estimators, due to the
 fact that the ch.f.\ admits a simple form. For example, $\phi_{\lambda}(\pi)=1/\cosh(\lambda\pi)=\beta$, say,  equals to the difference $\Pr(X_\lambda \ \mbox{even})-\Pr(X_\lambda \ \mbox{odd})$. This
 can be consistently and unbiasedly estimated by $\widehat{\beta}_n=n^{-1}\sum_{i=1}^n (-1)^{X_{i}}$, and
 a trivial application of the CLT leads to
 $\sqrt{n}(\widehat{\beta}_n-\beta)\stackrel{d}{\to}N(0,1-\beta^2)$, while the SLLN
 shows that $\widehat{\beta}_n$ is eventually positive w.p.\ $1$.
 Applying the delta-method (see \cite{12}) with $g(\beta)=\pi^{-1}\left(\log(1 + \sqrt{1 - \beta^2})-\log(\beta)\right)$, so that $g(\beta)=\lambda$, we obtain
 \[
 \sqrt{n}\left(g(\widehat{\beta}_n)-\lambda\right)\stackrel{d}{\to}N\Big(0,\cosh(\pi\lambda)^2/\pi^2\Big).
 \]
 However, compared to the MLE, the closed-form estimator $g(\widehat{\beta}_n)$ is by far less efficient.
 Thus, it is natural to seek for closed-form highly efficient estimators, and this may be
 possible as in the continuous case.
 In the continuous case it is shown that the asymptotic relative efficiency
 of the geometric mean of the absolute values of the observations is $8/\pi^2\simeq 81\%$, and
 in \cite{KP-N5} a more efficient closed-form estimate is proposed.
 Also, highly efficient estimators
 that are based on the ch.f.\  may be obtained by adapting the methodology of \cite{K-N4}
 to the present discrete case. However, such results are beyond the scope of the present note.
 Note that the Fisher information in the continuous Cauchy scale family
 equals to $1/(2\lambda^2)$ (compare to (\ref{I})), and the likelihood
 equation is as in (\ref{Like}),
 with the absence of the term $\pi\lambda/\sinh(2\pi\lambda)$.

  \begin{REM}
 \label{rem.2}
 The series in Corollary \ref{cor.2} is
 of some interest in itself,
 because of the computation of the sum
 $
 \sum_{n=1}^{\infty} (\lambda^2+n^2)^{-1}
 $
 in a closed form. Then, e.g., taking limits
 as $\lambda\searrow 0$,
 we arrive at the famous Euler sum, $\sum_{n=1}^\infty n^{-2}=\pi^2/6$.
 Moreover,
 differentiating term by term with respect to $\lambda$
 we can evaluate the series
 \[
 \sum_{n=1}^{\infty} \frac{1}{(\lambda^2+n^2)^2}.
 \]
 From this, taking limits
 as $\lambda\searrow 0$,
 we arrive at the sum for $\zeta(4)$,
 that is, $\sum_{n=1}^\infty n^{-4}=\pi^4/90$;
 clearly, this process can be continued to evaluate all
 $\zeta(2s)$ values, as well as the series
 $\sum_{n=1}^{\infty} (\lambda^2+n^2)^{-s}$, $s=1,2,\ldots$ .
 \end{REM}

 \begin{REM}
 \label{rem.3}
 Differentiating $m$ times with respect to $\lambda^2$
 the series in Lemma
 \ref{lem.2}, it is
 possible to introduce
 and investigate discrete Student-type
 families with $\nu=2m+1$ degrees of freedom, that is,
 p.m.f.'s
 of the form
 \be
 \label{5}
 f_{\nu;\lambda}(k)= \frac{c_{\nu;\lambda}}{(\lambda^2+k^2)^{(\nu+1)/2}},
 \ \ \ k\in\Z, \
 \ \ \nu=1,3,5,\ldots, \ \ \
 \lambda>0,
 \ee
 admitting
 closed-form
 ch.f.'s $\phi_{\nu;\lambda}(t)$ and
 explicit
 normalizing
 constants $c_{\nu;\lambda}$. However, the situation becomes
 quite complicated for even values of $\nu$.
 \end{REM}


 \section{Some distributional properties of the \mbox{${\cal CC}$} family}
 \label{sec.4}

 \noindent
 We observe that the ch.f.\ $\phi_\lambda(t)$ is not differentiable at the
 points $t=2k\pi$, $k\in\Z$ (c.f.\ Fig.\ \ref{fig.1}). It is known that
 a random variable $Y_1$ satisfies a weak law of large numbers,
 that is,
 \[
 \overline{Y}_n :=\frac{Y_1+\cdots+Y_n}{n}\to \mbox{ some constant } c,
 \mbox{
 in probability},
 \]
 if and only if its ch.f., $\phi_{Y_1}$, is differentiable at $t=0$; then,
 $\phi_{Y_1}'(0)=i c$ where $i$ is the imaginary unit (the problem
 was treated by A.\ Zygmund and E.J.G.\ Pitman, and it is closely connected to
 Khintchine's weak law of large numbers; see Feller 1966, p.\ 528
 and van der Vaart 1998, p.\ 15).
 Hence, the distributions of the \mbox{${\cal CC}$} family do not
 satisfy the weak law of large numbers, since their ch.f.'s are
 not differentiable at $t=0$. Therefore, it is of some interest
 to study the asymptotic behavior of the sample means from a \mbox{${\cal CC}$}
 random variable with p.m.f.\
 as in (\ref{3}).
 Recall the well-known continuous counterpart, which says that
 $\overline{X}_n$ is the same Cauchy for all $n$ (Cauchy r.v.'s are stable).

 We have the following result.
 \begin{theo}
 \label{theo.2}
 If $X_1,X_2,\ldots$ are independent identically distributed random variables
 with p.m.f.\ as in {\rm (\ref{3})} then
 \[
  \overline{X}_n\stackrel{d}{\to} \lambda\tanh(\lambda\pi) \ Z,
 \]
 where $Z$ is standard {\rm (continuous)} Cauchy with density
 {\rm (\ref{1b})}.
 \end{theo}
 \begin{pr}{Proof}
 Fix $t\geq0$.
 Theorem \ref{theo.3} shows that
 the ch.f.\ of $\overline{X}_n$ is
 given by
 \[
 \phi_\lambda(t/n)^n=
 \left(\frac{\cosh\big(\lambda
 (\pi-t/n)\big)}{\cosh(\lambda\pi)}\right)^n, \  n\geq \frac{t}{2\pi}.
 \]
 Using this, it is easy to verily (e.g., by taking logarithms)
 that $\phi_\lambda(t/n)^n\to e^{- c t}$,
 $t\geq 0$, where $c=\lambda\tanh(\lambda\pi)$.
 %
 Finally, from the fact that $\phi_\lambda$ is even, it follows
 that $\phi_\lambda(t/n)^n\to e^{-c|t|}$ for all
 $t\in\R$, which is the ch.f.\ of $cZ$, and the result follows from
 the continuity theorem of characteristic functions.
 \bigskip
 $\Box$
 \end{pr}

 Unlike the usual Cauchy scale family,
 the \mbox{${\cal CC}$} family is not
 convolution closed; however, it is ``almost"
 closed. More precisely, the following result
 holds.
 \begin{theo}
 \label{theo.4}
 For independent r.v.'s
 $X$, $Y$ in $\mbox{\rm \mbox{${\cal CC}$}}$ with $X\sim f_{\lambda_1}$ and
 $Y\sim f_{\lambda_2}$,
 the ch.f.\ of $X+Y$ is given by
 \[
 \phi_{X+Y}(t)=
 \frac{\alpha(\lambda_1+\lambda_2)}{2\alpha(\lambda_1)
 \alpha(\lambda_2)} \phi_{\lambda_1+\lambda_2}(t)
 +\frac{\alpha(|\lambda_2-\lambda_1|)}{2\alpha(\lambda_1)
 \alpha(\lambda_2)} \phi_{|\lambda_2-\lambda_1|}(t), \ \ t\in\R,
 \]
 where $\phi_0(t)\equiv 1$ is the ch.f.\ of the degenerate
 r.v.\ $X_0$ with $\Pr(X_0=0)=1$, and
 $\alpha(\lambda):=\cosh(\lambda\pi)$, $\lambda\geq0$.
 Consequently, $X+Y$ is a mixture of
 two r.v.'s that are members of {\rm \mbox{${\cal CC}$}} family,
 \[
 \Pr(X+Y=k)=\frac{\alpha(\lambda_1+\lambda_2)}{2\alpha(\lambda_1)
 \alpha(\lambda_2)} f_{\lambda_1+\lambda_2}(k)
 +\frac{\alpha(|\lambda_2-\lambda_1|)}{2\alpha(\lambda_1)
 \alpha(\lambda_2)} f_{|\lambda_2-\lambda_1|}(k), \ \ k\in\Z.
 \]
 \end{theo}
 \begin{pr}{Proof}
 Set
 \[
 p=\frac{\alpha(\lambda_1+\lambda_2)}{2\alpha(\lambda_1)
 \alpha(\lambda_2)}, \ \ \ \ q=\frac{\alpha(|\lambda_2-\lambda_1|)}{2\alpha(\lambda_1)
 \alpha(\lambda_2)}.
 \]
 Obviously, $p>0$ and $q>0$. Also, using the formula
 \be
 \label{4}
  \cosh(x)\cosh(y)=\frac{1}{2}\cosh(x+y)+\frac{1}{2}\cosh(y-x)
 \ee
 it is easily seen that  $p+q=1$. Restricting our attention to the interval
 $0\leq t\leq 2\pi$, we have
 \[
  \phi_{X+Y}(t)=\phi_{\lambda_1}(t)\phi_{\lambda_2}(t)
  =\frac{\cosh(\lambda_1(t-\pi))\cosh(\lambda_2(t-\pi))}{\alpha(\lambda_1)
  \alpha(\lambda_2)}
 \]
 and a final application of (\ref{4}) to the numerator, taking into
 account Theorem \ref{theo.3}, completes the proof.
 $\Box$
 \bigskip
 \end{pr}

 \begin{REM}
 \label{rem.1}
 If $X,Y$ are i.i.d.\ from $f_\lambda$ then, since $\alpha(0)=1$
 and $f_0(k)=I(k=0)$, we get
 \[
 \Pr(X+Y=k)=\left\{
 \begin{array}{ll}
 \displaystyle
 \frac{1}{2\cosh(\lambda\pi)^2}+\frac{\tanh(\lambda\pi)}{2\lambda\pi},
  & k=0,
 \\
 &
 \\
 \displaystyle
 \frac{\tanh(\lambda\pi)}{\pi} \ \frac{2\lambda}{(2\lambda)^2+k^2},
 & k\in \Z^*.
 \end{array}
 \right.
 \]
 This formula quantifies the fact that the p.m.f.\ of $X+Y$ lies outside \mbox{${\cal CC}$},
 but it is close, in some sense, to $f_{2\lambda}$; in fact, the ratio $f_{X+Y}(k)/f_{2\lambda}(k)$
 does not vary with $k\in\Z^*$.
 \bigskip
 \end{REM}

 A ch.f.\ $\phi$ (or the corresponding r.v.\ $X$) is called infinitely divisible (i.d.)
 if for each $n$, we can find a ch.f.\ $\phi_n$ such that $\phi_n^n=\phi$;
 equivalently, if $X_{1,n}+\cdots+X_{n,n}$ has the same distribution as $X$, where
 $X_{1,n},\ldots,X_{n,n}$ are i.i.d.\ with ch.f.\ $\phi_n$.
 Properties of this kind are included in what is called
 "arithmetic of probability laws" (\cite{Lu-N6},
 \cite{Lu-N7}), and a vast
 bibliography exists, see, e.g.,
 \cite{GH-N2},
 \cite{Cr-N1},
 \cite{KS-N3},
 \cite{Lu-N6},
 \cite{Lu-N7},
 \cite{SvH-N8}, and references therein.

  Since the notion of i.d.\ is related to limit theorems  of sums
  of independent r.v.'s, it would be useful to know whether
  the \mbox{${\cal CC}$} family is i.d. This is indeed the case, and it follows immediately from a result
  of Polya, because the
  ch.f.\ $\phi_{\lambda}$ is even, log-convex in $[0,2\pi]$ and $2\pi$ periodic,
  see \cite{KS-N3}, \cite{Lu-N7}.
  In fact, $\phi_\lambda^\alpha$ is a ch.f.\ for all $\lambda\geq0$ and $\alpha\geq0$.

  As is well known, the notion of self-decomposability, as well as that of stability, do not apply to discrete r.v.'s. Recall that $X$ is stable if, for each $n$, we can find constants $\alpha_n>0$ and $\beta_n\in\R$ such that $X$ and $(X_1+\cdots+X_n)/\alpha_n-\beta_n$ have the same distribution,
  where $X_1,\ldots,X_n$ are i.i.d.\ copies $X$. Obviously, the class of stable distributions
  is a proper subset of i.d.\ distributions. Due to a fundamental result of L\'{e}vy,
  stable distributions are very important because their class
  contains exactly  all possible limits of (properly) normalized sums of i.i.d. r.v.'s.
  Every stable distribution has a ch.f.\ that can be expressed in a closed form, and the
  corresponding r.v.\ is absolutely continuous. The subclass of symmetric stable ch.f.'s,
  after a location-scale transformation, can be written as ${\cal S}=\{\phi_\alpha(t)=e^{-|t|^\alpha}$,  $0<\alpha\leq2\}$. Only the densities that correspond to $\alpha=1/2$ (L\'{e}vy),
  $\alpha=1$ (Cauchy) and $\alpha=2$ (Normal), have known explicit forms.

  It is natural to ask whether the \mbox{${\cal CC}$} family contains discrete stable distributions,
  in the sense
  of \cite{SvH-N8}. However, the definitions in \cite{SvH-N8} are designed for non-negative
  integer-valued r.v.s, and are based on probability generating functions; it is not obvious
  how to extend these results to the \mbox{${\cal CC}$} case. The following definition
  provides a different approach that seems to be natural for our case.

  \begin{DEFI}
  \label{def.2}
  Let $\Lambda$
  be a set of indices, consider a parametric family
  ${\cal F}=\{\phi_\lambda, \ \lambda\in \Lambda\}$  of discrete, integer-valued, ch.f.'s, and
  let ${\cal F}'$ be the corresponding family of random variables.
  Then, ${\cal F}$
  is called discrete stable (DSF)
  if for each $\phi_\lambda\in {\cal F}$, we can find a sequence
  of indices $\{\lambda_n\}_{n=1}^{\infty}\subset \Lambda$
  such that $\phi_{\lambda_n}^n\to\phi_{\lambda}$.
  Equivalently, if every random variable in ${\cal F}'$
  is the weak limit of sums of i.i.d.\ r.v.'s from ${\cal F}'$.
  \end{DEFI}

  The usual Poisson family is DSF, as well as the Negative Binomial.
  In order for
  such a model to be useful in practice, the family ${\cal F}$ should not contain "too many" ch.f.'s.
  Also, it is plausible to consider those DSF's that satisfy some kind of discrete attraction, in the
  sense that (non-normalized) sums of several i.i.d.\ discrete r.v.'s converge weakly to one of the members of the DSF.
  It is clear that the Compound Poisson that is produced by a fixed discrete ch.f.\ $\psi$, namely,
  ${\cal F}=\{\phi_{\lambda}(t)=e^{\lambda(\psi(t)-1)}, \ \lambda\geq 0\}$, is such a useful DSF model.
  On the other hand, the complete Compound Poisson model (allowing any $\psi$ in the exponent) seems to
  be too wide.
  Regarding the \mbox{${\cal CC}$} family we have the following result.

  \begin{theo}
  \label{theo.6}
  The \mbox{${\cal CC}$} family is not DSF. To be more specific, suppose
  $\{\phi_{\lambda_n}\}_{n=1}^{\infty}\subset \mbox{\mbox{${\cal CC}$}}$
  where $\lambda_n\geq 0$
  is an arbitrary sequence, and $\phi_{\lambda_n}$ is as in Theorem {\rm \ref{theo.3}}.
  Then, {\rm (i)} and {\rm (ii)} below are equivalent.
  \vspace{1ex}

  \noindent
  {\rm (i)} There is a point $t_0\in(0,2\pi)$
  such that $\lim_n \phi_{\lambda_n}(t_0)^n=\delta>0$.
  \vspace{1ex}

  \noindent
  {\rm (ii)} It holds $\lambda_n=\theta/\sqrt{n}+o(1/\sqrt{n})$, where $\theta=(-2\log\delta)^{1/2}(t_0(2\pi-t_0))^{-1/2}\geq 0$.
  \vspace{1ex}

  \noindent
  If {\rm (i)} or {\rm (ii)} is satisfied then
  $\phi_{\lambda_n}(t)^n\to\psi(t):=\exp(-\theta^2 t(2\pi-t)/2)$ uniformly
  in $t$, $0\leq t\leq 2\pi$, and the limiting ch.f.\
  $\psi$
  {\rm (extended to be $2\pi$-periodic)} is an infinitely divisible ch.f.
  \end{theo}

  Before proving Theorem \ref{theo.6}, we provide some remarks.
  The limiting ch.f.\ $\psi$ is a
  Compound Poisson one. Indeed, the exponent can be written as
  $\lambda(\psi_1(t)-1)$, where $\psi_1(t)=1-\theta^2 t(\pi-t/2)/\lambda$
  and, e.g., $\lambda\geq\pi^2\theta^2/2$ (we shall see below that
  the minimum value of $\lambda$ for which $\psi_1$ is a ch.f.\ is $\lambda_0=\pi^2\theta^2/3$).
  Then, it follows that the even, $2\pi$-periodic
  function $\psi_1$ is nonnegative, decreasing and convex in $[0,\pi]$, and so, by Polya's
  sufficiency criterion
  (see \cite{KS-N3}) it is a ch.f.\ of an integer-valued r.v.
  Clearly, the parametric family produced by all possible limits
  from \mbox{${\cal CC}$}, namely,
  ${\cal F}=\{\psi_{\lambda}(t)=e^{-\lambda t(2\pi-t)}, \ \lambda\geq0, \ t\in[0,2\pi]\}$,
  forms a DSF according to Definition \ref{def.2}. By applying the inversion formula
  for ch.f.'s of integer-valued r.v.'s,
  namely,
  \[
  \Pr(X=k)=\frac{1}{2\pi}\int_{-\pi}^{\pi} e^{-ikt}
  \phi_X(t)dt, \ \ k\in\Z,
  \]
  it is recognized that the p.m.f.'s
  in ${\cal F}$
  do not admit closed forms. Indeed, if $Y_\lambda\sim\psi_{\lambda}$
  then the preceding formula reduces to
  \[
  \Pr(Y_\lambda=k)=\frac{1}{\pi}\int_{0}^{\pi} \cos(kt)
  e^{-\lambda t(2\pi-t)}dt, \ \ k\in\Z,
  \]
  and this integral cannot be computed in terms of elementary functions (unless $\lambda=0$).
  Moreover, if we make use of the preceding formula with $\psi_1$ instead of $\psi$,
  we can easily obtain the p.m.f.\ of the r.v.\ $W$ with ch.f.\ $\psi_1$.
  Setting for convenience $c=\theta^2/\lambda$
  one finds $\Pr(W=0)=1-c\pi^2/3$ (so that $c\leq 3/\pi^2$ and, hence, $\lambda\geq \theta^2\pi^2/3$)
  and $\Pr(W=k)=c/k^2$, $k\in \Z^*$. According to Theorem \ref{theo.6},
  these remarks provide
  a detailed description of the class of the limiting distributions of sums of i.i.d.\
  r.v.'s from \mbox{${\cal CC}$}.

  The following lemma will be used in the proof of Theorem \ref{theo.6}.
  \begin{lem}

  \label{lem.3}
  {\rm (i)}
  Let $\{\beta_n\}_{n=1}^\infty\subset (0,1]$, assume that $\beta_n^n\to\beta\in(0,1]$
  and set $B=-\log \beta$.
  Then,
  $\beta_n=1-B/n+o(1/n)$.
  \medskip

  \noindent
  {\rm (ii)}
  Fix $x_0\in[0,1)$, and define the function $f(y):=\cosh(x_0 y )/\cosh(y)$, $y\geq 0$.
  Suppose that $\{\alpha_n\}_{n=1}^\infty\subset [0,\infty)$
  and that $f(\alpha_n)^n\to\delta\in(0,1]$. Then, $\alpha_n=\alpha/\sqrt{n}+o(1/\sqrt{n})$, where
  $\alpha=\sqrt{(-2\log\delta)/(1-x_0^2)}$. 
  \end{lem}
  \begin{pr}{Proof}
  (i) Despite the fact that (i) is known, we provide a very quick proof here.
  The inequality $y\leq -\log(1-y)\leq y/(1-y)$ ($0\leq y<1$), applied $y=1-\beta_n$, yields
  $
  \beta_n(-n \log\beta_n)\leq n(1-\beta_n)\leq -n\log\beta_n,
  $
  and since the upper bound implies that $\beta_n\to 1$, both bounds
  converge to $B$.
  %
  \\
  (ii) The sequence $n\alpha_n^2$ is bounded. Indeed, assuming the contrary,
  it follows that for any $M>0$ (arbitrarily large) we can find a subsequence $n_k$
  such that $\alpha_{n_{k}}>M/\sqrt{n_{k}}$ for all $k$.
  Since it is easily checked that $f'(y)<0$ for $y>0$, the positive continuous
  function $f$ is strictly decreasing, with $f(0)=1$, $f(\infty)=0$ (recall that $0\leq x_0<1$).
  Therefore, $f\big(\alpha_{n_k}\big)^{n_k}\leq f\big(M/\sqrt{n_{k}}\big)^{n_k}\to\exp\big(-M^2(1-x_0)^2/2\big)$,
  as $k\to\infty$. Thus, $\liminf f(\alpha_n)^n\leq \exp\big(-M^2(1-x_0)^2/2\big)$, and since
  $M>0$ is arbitrary, $\liminf f(\alpha_n)^n\to 0$. This contradicts the hypothesis
  $f(\alpha_n)^n\to \delta>0$, and verifies that the sequence $n\alpha_n^2$ is, indeed, bounded.
  Hence, $\alpha_n\to 0$. By applying a Taylor development to the function $f$ it can be checked
  that for
  $y\geq 0$,  sufficiently close to zero,
  \[
  1-\frac{1}{2}(1-x_0^2)y^2\leq  f(y) \leq 1-\frac{1}{2}(1-x_0^2)
  y^2+\frac{1}{24}(1-x_0^2)(5-x_0^2)y^4, \ \ 0\leq y<\epsilon.
  \]
  Substituting $y=\alpha_n$ (which tends to zero) we obtain the inequality
  \[
  A n(1-f(\alpha_n))\leq n\alpha_n^2\leq A n(1-f(\alpha_n))+B \alpha_n^2(n\alpha_n^2),
  \ \ n\geq n_0,
  \]
  with $A=2/(1-x_0^2)$, $B=(5-x_0^2)/12$. Since $f(\alpha_n)^n\to\delta\in(0,1]$
  (and $0<f(\alpha_n)\leq 1$), it follows from part (i)
  that $n(1-f(\alpha_n))\to-\log\delta$, and the preceding inequality
  shows that $n\alpha_n^2\to(-\log\delta)A$, completing the proof.
  \bigskip
  $\Box$
  \end{pr}

  \noindent
  \begin{pr}{Proof of Theorem \ref{theo.6}}
  Assume first that (ii) holds, that is,
  $\lambda_n=\theta/\sqrt{n}+o(1/\sqrt{n})$ for some $\theta\geq 0$.
  It is straightforward to verify
  that $\phi_{\lambda_n}(t)^n$ converges pointwise to $\psi(t)$ as given, and from the
  fact that $\psi$ is continuous at the origin, the convergence is uniform
  at compacts, and in particular, in $[0,2\pi]$. Obviously, (i) is satisfied for (any
  choice of) $t_0\in(0,2\pi)$
  with $\delta=\psi(t_0)=\exp(-\theta^2 t_0(2\pi-t_0)/2)>0$.

  Assume now that (i) holds, i.e., suppose that for
  a fixed $t_0\in(0,2\pi)$,  $\phi_{\lambda_n}(t_0)^n\to\delta>0$.
  Due to symmetry ($\phi_{\lambda_n}(t)=\phi_{\lambda_n}(2\pi-t)$),
  we can further assume that $0<t_0\leq \pi$. Set $\alpha_n=\pi\lambda_n$,
  $x_0=1-t_0/\pi\in[0,1)$,  and consider
  the function $f(y)=\cosh(x_0y)/\cosh(y)$, $y\geq 0$, as in Lemma \ref{lem.3}.
  Then, $\phi_{\lambda_n}(t_0)=f(\alpha_n)$, and
  by assumption, $f(\alpha_n)^n\to\delta>0$ (certainly, $\delta\leq 1$).
  Hence, from Lemma \ref{lem.3}(ii) we conclude that
  $n\alpha_n^2\to (-2\log\delta)/(1-x_0^2)$, that is,
  $n\lambda_n^2\to (-2\log\delta)/(t_0(2\pi-t_0))$, which verifies (ii).
  \bigskip
   $\Box$
  \end{pr}

  It is of some interest to observe that, according to Theorem \ref{theo.6},
  the limiting ch.f.\ exists
  if we can merely show
  the convergence $\phi_{\lambda_n}(t_0)^n\to\delta>0$
  for a single nontrivial point $t_0$ (i.e., $t_0\neq 2k\pi$).
  Then, $\psi(t)$ is uniquely determined from
  the pair $(t_0,\delta)$,  Also,
  the limiting distribution is degenerate at zero if and only if $\delta=1$
  (which is corresponds to $\theta=0$ in Theorem \ref{theo.6}(ii)).

  Another related problem concerns the extended \mbox{${\cal CC}$} class,
  defined as the family of ch.f.'s
  $\mbox{\mbox{${\cal CC}$}}^+:=\{\phi^\alpha:\ \phi\in\mbox{\mbox{${\cal CC}$}},\ \alpha>0 \}$.
  Since every $\phi\in\mbox{\mbox{${\cal CC}$}}$ is $2\pi$-periodic, decreases in $[0,\pi]$
  and is log-convex in $[0,2\pi]$, the same is true for
  all ch.f.'s in $\mbox{\mbox{${\cal CC}$}}^+$. Hence,
  $\mbox{\mbox{${\cal CC}$}}^+$ is a family of i.d.\
  ch.f.'s. This family is similar to the (continuous) Cauchy scale family.
  Cram\'{e}r  \cite{Cr-N1} showed that all stable centered distributions
  with exponent $\alpha<2$ are not factor closed. This means that, e.g.,
  the ch.f.\ of the standard Cauchy, $e^{-|t|}$,
  can be written as $\phi_1\phi_2$, with $\phi_i$ ($i=1,2$) lying outside the class of
  Cauchy ch.f.'s. So, it is fairly expected that the same is true
  for $\mbox{\mbox{${\cal CC}$}}^{+}$. Indeed, it can be proved that this is the case, and, as a concrete
  example, we provide the following $2\pi$-periodic $\phi_i$'s:
  \begin{eqnarray*}
  \phi_1(t)&=&\left(\frac{\cosh(t-\pi)}{\cosh(\pi)}\right)^{1/2}
  \left(\frac{1+\pi^4}{1+(t-\pi)^4}\right)^{1/50}, \ \   0\leq t\leq 2\pi,
  \\
  \phi_2(t)&=&\left(\frac{\cosh(t-\pi)}{\cosh(\pi)}\right)^{1/2}
  \left(\frac{1+(t-\pi)^4}{1+\pi^4}\right)^{1/50},
  \ \   0\leq t\leq 2\pi.
  \end{eqnarray*}
  It can be checked that both functions are positive, decreasing in $[0,\pi]$,
  and convex ($\phi_1$ is log-convex) in $[0,2\pi]$ and hence,
  their $2\pi$-periodic extensions (which are even functions)
  are ch.f.'s, see \cite{Lu-N7}.
  Obviously, these ch.f.'s lie outside $\mbox{\mbox{${\cal CC}$}}^+$,
  and, trivially, their product equals to the standard discrete Cauchy ch.f.\ of Theorem \ref{theo.1}.
  \bigskip

  {\bf Acknowledgements.} I would like to thank an anonymous referee who
  motivated several of the results that are presented in Section \ref{sec.4},
  regarding infinite divisibility and stability in the discrete case.



\begin{thebibliography}{99}
 {\small
 \bibliographystyle{Chicago}

 \bibitem{1}
 Arnold, B.C.; Brockett, P.L.\ (1992).
 On distributions whose component ratios are Cauchy.
 {\it Amer.\ Statist.}, {\bf 46}, 25--26.
 \vspace{-1ex}




 \bibitem{4}
 Cacoullos, T.\ (2014).
 Polar angle tangent vectors follow Cauchy distributions under
 spherical symmetry. {\it J.\ Multivariate Anal.}, {\bf 128}, 147--153.
 \vspace{-1ex}



 \bibitem{Cr-N1}
 Cram\'{e}r, H.\ (1949). On the factorization of certain probability distributions,
 {\it Aekiv f\:{o}r Matematik}, {\bf 7}, 61--65.
 \vspace{-1ex}

 \bibitem{6}
 Feller, W.\ (1966). {\it An Introduction to Probability Theory and its
 Applications}, Vol.\ II. Wiley, N.Y.
 \vspace{-1ex}

 \bibitem{GH-N2}
 Geluk, J.L; de Haan, L.\ (2000).
 Stable probability distributions and their
 domains of attraction: a direct approach
 {\it Probab. Math. Statist.},
 {\bf 20}, 169--188.
 \vspace{-1ex}

 \bibitem{8}
 Jones,  M.C.\ (1999). Distributional relationships arising from
 simple trigonometric formulas. {\it Amer.\ Statist.},
 {\bf 53}, 99--102.
 \vspace{-1ex}


 \bibitem{10}
 Jones, M.C.\ (2008). The distribution of the ratio $X/Y$ for all centred
 elliptically symmetric distributions.
 {\it J.\ Multivariate Anal.}, {\bf 99}, 572--573.
 \vspace{-1ex}

 \bibitem{KS-N3}
 Keilson, J.; Steutel, F.W. (1972).
 Families of infinitely divisible distributions closed under mixing and convolution.
 {\it Ann.\ Math.\ Statist.}, {\bf 43}, 242-250.
 \vspace{-1ex}

 \bibitem{K-N4}
 Koutrouvelis. I.A.\ (1982). Estimation of location and scale in Cauchy distributions using the empirical characteristic function.
 {\it Biometrika}, {\bf 69}, 205--213.
 \vspace{-1ex}

 \bibitem{KP-N5}
 Kravchuk, O.Y; Pollett, P.K.\ (2012).
 Hodges-Lehmann scale estimator for Cauchy
 distribution.
 {\it Commun.\ Statist.--Theory Meth.}, {\bf 41}, 3621--3632.
 \vspace{-1ex}

 \bibitem{LG-N9} Lehmann, E.L.; Gasella, G.\ (1998). {\it Theory of Point
 Estimation}, 2nd ed. Springer.,
  \vspace{-1ex}
 N.Y.

 \bibitem{Lu-N6}
 Lukacs, E.\ (1961). Recent developments in the theory of characteristic functions.
 {\it Proc.\ Fourth Berkeley Symposium}, {\bf 2}, 307--335.
 \vspace{-1ex}

 \bibitem{Lu-N7}
 Lukacs, E.\ (1972).
 A survey of the theory of characteristic functions.
 {\it Adv.\ Appl.\ Probab.}, {\bf 4}, 1--38.
 \vspace{-1ex}




 \bibitem{11}
 Renshaw, E.\ (2011).
 {\it Stochastic Population Processes: Analysis, Approximations, Simulations.}
 Oxford University Press.
 \vspace{-1ex}

 \bibitem{SvH-N8}
 Steutel, F.W.; van Harn, K.\ (1979).
 Discrete analogues of self-decomposability and stability.
 {\it Ann.\ Probab.}, {\bf 7}, 893--899.
 \vspace{-1ex}


 \bibitem{12}
 van der Vaart, A.W.\ (1998). {\it Asymptotic Statistics}. Cambridge
 University Press, N.Y.
 \vspace{-1ex}



 }
 \end{thebibliography}
 \end{document}